\documentclass[]{amsart}
\usepackage{amssymb}

\theoremstyle{plain}
\newtheorem{theorem}{Theorem}

\newtheorem{lemma}[theorem]{Lemma}
\newtheorem{proposition}[theorem]{Proposition}

\theoremstyle{definition}

\newtheorem{remark}[theorem]{Remark}

\newcommand{\abs}[1]{\lvert#1\rvert}
\newcommand{\norm}[1]{\lVert#1\rVert}
\newcommand{\bigabs}[1]{\bigl\lvert#1\bigr\rvert}

\renewcommand{\le}{\leqslant}
\renewcommand{\ge}{\geqslant}

\newcommand{\term}[1]{{\textit{\textbf{#1}}}}

\def\Orb{{\rm Orb}}
\def\Range{{\rm Range}}

\begin{document}

\title[Spectrum of w-hypercyclic operator]
      {Spectrum of a weakly hypercyclic\\
       operator meets the unit circle}

\author[S.J.~Dilworth, V.G.~Troitsky]{S. J. Dilworth and Vladimir G.~Troitsky} 
\address{Department of Mathematics, University of South Carolina,
         Columbia, SC 29208. USA.}
\curraddr{Department of Mathematics, The University of Texas at Austin, Austin,
          Texas 78712. USA}
\address{Department of Mathematical and Statistical Sciences, 632 CAB,
  University of Alberta, Edmonton, AB T6G\,2G1. Canada.}
\email{dilworth@math.sc.edu\\vtroitsky@math.ualberta.ca}

\thanks{The research of the first author was supported by a summer
  research grant from the University of South Carolina while on
  sabbatical as a Visiting Scholar at The University of Texas at
  Austin.}

\keywords{Weakly hypercyclic operator, spectrum, spectral set, weak closure.}
\subjclass{Primary: 47A16,47A10,47A25}

\begin{abstract}
  It is shown that every component of the spectrum of a weakly
  hypercyclic operator meets the unit circle. The proof is based on
  the lemma that a sequence of vectors in a Banach space whose norms
  grow at geometrical rate doesn't have zero in its weak closure.
\end{abstract}

\maketitle



Suppose that $T$ is a bounded operator on a nonzero Banach space $X$.
Given a vector $x\in X$, we say that $x$ is \term{hypercyclic} for $T$
if the orbit $\Orb_Tx=\{T^nx\}_n$ is dense in $X$. Similarly, $x$ is
said to be \term{weakly hypercyclic} if $\Orb_Tx$ is weakly dense in
$X$. A bounded operator is called \term{hypercyclic} or \term{weakly
  hypercyclic} if it has a hypercyclic or, respectively, a weakly
hypercyclic vector. It is shown in~\cite{Chan:02} that a weakly
hypercyclic vector need not be hypercyclic, and there exist weakly
hypercyclic operators which are not hypercyclic.  C.~Kitai
showed in~\cite{Kitai:82} that every component of the spectrum of a
hypercyclic operator intersects the unit circle. K.~Chan and
R.~Sanders asked in~\cite{Chan:02} if the spectrum of a weakly
hypercyclic operator meets the unit circle.  In this note we show that
every component of the spectrum of a weakly hypercyclic operator meets
the unit circle.

\begin{lemma} \label{l:fast}
  Let $X$ be a Banach space and let $c>1$. Suppose that $x_n \in X$
  satisfies $\norm{x_n}\ge c^n$ for all $n\ge1$. Then 
  $0\notin\overline{\{x_n\}_n}^w$.
\end{lemma}

\begin{proof}
  Let $N$ be the smallest positive integer such that $c^N>2$.  We
  shall prove that there exist $F_1,\dots,F_N \in X^*$ such that
  \begin{equation} \label{eq: F1toFn}
    \max_{1\le k \le N}\bigabs{F_k(x_n)} \ge 1 \qquad (n\ge1). 
  \end{equation}
  Since $\norm{x_n} \ge c^n$, by replacing $x_n$ by $(c^n/\norm{x_n}) x_n$,
  it suffices to prove \eqref{eq: F1toFn} for the case in which
  $\norm{x_n} = c^n$ for all $n\ge1$. First suppose that $c>2$, so that
  $N=1$.  We have to construct $F_1 \in X^*$ such that $\bigabs{F_1(x_n)} \ge
  1$ for all $n\ge1$.  First choose $f_1 \in X^*$ with $f_1(x_1)=1$.
  Then either $\bigabs{f_1(x_2)}<1$ or $\bigabs{f_1(x_2)} \ge 1$.  In the former
  case the Hahn-Banach theorem guarantees the existence of $g_2 \in
  X^*$ such that $\norm{g_2} \le 1/\norm{x_2} = c^{-2}$, $\bigabs{g_2(x_2)} = 1-
  \bigabs{f_1(x_2)}$, and $\bigabs{(f_1+g_2)(x_2)} =1$. In the latter case, set
  $g_2=0$. Note that
  $$\bigabs{(f_1+g_2)(x_1)} \ge 1 - \norm{g_2}\norm{x_1} \ge 1- c^{-1}.$$
  Set $f_2=f_1+g_2$. Repeating this argument, we can find $g_3 \in X^*$ such that
  $\norm{g_3} \le 1/\norm{x_3} = c^{-3}$ and $\bigabs{(f_2+g_3)(x_3)}\ge1$. Note that
  $$\bigabs{(f_2+g_3)(x_1)}\ge \bigabs{f_2(x_1)} - \norm{g_3}\norm{x_1}
    \ge 1 - c^{-1} - c^{-2}$$
  and also that
  $$\bigabs{(f_2+g_3)(x_2)} \ge 1 - \norm{g_3}\norm{x_2} \ge 1 - c^{-1}.$$
  Set
  $f_3=f_2+g_3$. Continuing in this way we obtain $f_n\in X^*$ 
  such that (setting $g_n = f_n-f_{n-1}$)  $\norm{g_n} \le c^{-n}$ and
  \begin{equation} \label{eq: induct}
    \bigabs{f_n(x_k)} \ge 1 - \sum_{i=1}^{n-k} c^{-i} \qquad(1 \le k \le n). 
  \end{equation}
  Thus, $\{f_n\}_n$ is norm-convergent in $X^*$ to some $f \in X^*$.
  From \eqref{eq: induct}, we obtain (since $c>2$)
  $$\bigabs{f(x_k)} = \lim_n \bigabs{f_n(x_k)} \ge 1 - \sum_{i=1}^\infty
      c^{-i}= \frac{c-2}{c-1} >0.$$
  Set $F_1 = (c-1)(c-2)^{-1}f$, to complete the proof in the case $c>2$.
  
  Now suppose that $1< c <2$. Set $\alpha = c^N >2$. For each $1
  \le k \le N$, consider the sequence $y_n = x_{k+(n-1)N}$ ($n\ge1$).
  Then $\norm{y_n} = (c^k/\alpha) \alpha^n$.  Since $\alpha >2$ there
  exists $F_k \in X^*$ such that $\bigabs{F_k(y_n)}\ge 1$ for all $n \ge 1$,
  which proves \eqref{eq: F1toFn}.
\end{proof} 

\begin{remark} \label{rem: norming}
  Recall that a closed subspace $Y$ of $X^*$ is said to be
  \term{norming} if there exists $C>0$ such that
  $$\norm{x} \le C\sup\bigl\{\bigabs{f(x)}\:\colon\: f\in
  Y,\,\norm{f}\le1\bigr\}\qquad (x \in X).$$
  The argument of
  Lemma~\ref{l:fast} easily generalizes to give the following result.
  Suppose that $Y$ is norming for $X$ and that $\{x_n\}_n$ is a sequence
  in $X$ satisfying $\norm{x_n}\ge c^n$, where $c>1$. Then $0$ does
   not belong to the $\sigma(X,Y)$-closure of  
  $\{ x_n\}_n$. In
  particular, Lemma~\ref{l:fast} is valid for the weak-star topology
  when $X$ is a dual space.
\end{remark} 

We also make use of the following simple numerical fact. If $(t_n)$ is a
sequence in ${\mathbb R}^+\cup\{\infty\}$, then
\[
  \limsup\limits_{n\to\infty}\sqrt[n]{t_n}=
  \inf\{\,\nu>0 \mid
    \lim\limits_{n\to\infty}\tfrac{t_n}{\nu^n}=0\,\}=
  \inf\{\,\nu>0 \mid
    \limsup\limits_{n\to\infty}\tfrac{t_n}{\nu^n}
    <\infty\,\}.
\]
In particular, if $T$ is a bounded operator with spectral radius
$r$, then the Gelfand formula
$\lim_n\sqrt[n]{\norm{T^n}}=r$ yields that $\frac{\norm{T^n}}{\lambda^n}\to 0$ 
for every scalar $\lambda$ with $\abs{\lambda}>r$.

\begin{theorem}
  If $T$ is weakly hypercyclic, then every component of $\sigma(T)$
  meets $\{z\,:\,\abs{z}=1\}$.
\end{theorem}

\begin{proof}
  Let $x$ be a weakly hypercyclic vector for $T$. Let $\sigma$ be a
  non-empty component of $\sigma(T)$, denote
  $\sigma'=\sigma(T)\setminus\sigma$. Denote by $X_\sigma$ and
  $X_{\sigma'}$ the corresponding spectral subspaces, then $X_\sigma$
  and $X_{\sigma'}$ are closed, $T$-invariant, and $X=X_\sigma\oplus
  X_{\sigma'}$. Also, $\sigma(T_{|X_\sigma})=\sigma$ and
  $\sigma(T_{|X_{\sigma'}})=\sigma'$. Note that $\sigma'$ might be
  empty, in which case we have $X_\sigma=X$ and $X_{\sigma'}=\{0\}$.
  
  Denote by $P_\sigma$ the spectral projection corresponding to $\sigma$,
  then $X_\sigma=\Range P_\sigma$. Denote $y=P_\sigma x$. Without loss
  of generality, $\norm{y}=1$. Since $P_\sigma$ is bounded and,
  therefore, weakly continuous, and
  $\Orb_Ty=P_\sigma\bigl(\Orb_Tx\bigr)$, we conclude that $\Orb_Ty$ is
  weakly dense in $X_\sigma$. Thus, $y$ is weakly hypercyclic for
  $T_{|X_\sigma}$.
  
  Observe that the inclusion $\sigma\subseteq\{z\,:\,\abs{z}<1\}$ is
  impossible.  Indeed, in this case the spectral radius of
  $T_{|X_\sigma}$ would be less than 1, so that $T^ny\to 0$, which
  contradicts $y$ being weakly hypercyclic for $T_{|X_\sigma}$.
  
  Finally, we show that the inclusion
  $\sigma\subseteq\{z\,:\,\abs{z}>1\}$ is equally impossible. In this
  case $0\notin\sigma=\sigma(T_{|X_\sigma})$, so that $T_{|X_\sigma}$
  is invertible. Denote the inverse by~$S$.  Then $S$ is a bounded
  operator on $X_\sigma$ and by the Spectral Mapping Theorem
  $$\sigma(S)=\bigl\{\lambda\mid\lambda^{-1}\in\sigma(T_{|X_\sigma})\bigr\}
  \subset\bigl\{z\,:\,\abs{z}<1\bigr\}.$$ Therefore, $r(S)<a$ for some
  $0<a<1$. This yields $\lim_n\frac{\norm{S^n}}{a^n}=0$, so that
  $\norm{S^n}\le a^n$ for all sufficiently large $n$. In particular,
  $$1=\norm{y}=\norm{S^nT^ny}\le a^n\norm{T^ny},$$ so that
  $\norm{T^ny}\ge\frac{1}{a^n}$. Lemma~\ref{l:fast} asserts that
  $0\notin\overline{\{T^ny\}_n}^w$, which contradicts
  $y$ being  weakly hypercyclic for $T_{|X_\sigma}$.
\end{proof}

\begin{proposition} 
  Suppose that $Y$ is norming for $X$. If $T$ has a hypercyclic vector
  for the $\sigma(X,Y)$ topology, then the spectrum of $T$ intersects
  $\{z\,\colon\abs{z}=1\}$.
\end{proposition}

\begin{proof}
  Suppose, to derive a contradiction, that $\sigma(T)$ does not
  intersect the unit circle. We use the notation introduced above with
  $\sigma=\sigma(T) \cap \{z\,\colon\abs{z}<1\}$ and $\sigma'=\sigma(T)
  \setminus \sigma$.  Let $x$ be a hypercyclic vector for the
  $\sigma(X,Y)$-topology. Then $x = y+z$, where $y \in X_\sigma$ and
  $z \in X_{\sigma'}$. Since $\norm{T^ny} \rightarrow 0$, it follows
  easily that $z$ is hypercyclic, and hence that $z \ne 0$.
  But then there exists $c>1$ such that
  $\norm{T^nz}\ge c^n$ for all sufficiently large $n$, which
  contradicts Remark~\ref{rem: norming}.
\end{proof}


\begin{thebibliography}{CS}

\bibitem[CS]{Chan:02}
Kit~C. Chan and Rebecca Sanders.
\newblock A weakly hypercyclic operator that is not norm hypercyclic.
\newblock To appear in J. Operator Theory.

\bibitem[K]{Kitai:82}
Carol Kitai.
\newblock {\em Invariant Closed Sets for Linear Operators}.
\newblock PhD thesis, Univ. of Toronto, 1982.

\end{thebibliography}

\end{document}